\documentclass[11pt,reqno]{amsart}
\usepackage{amssymb,graphicx}

\usepackage{amsfonts,amscd,amsthm,amsmath}

\usepackage{cite}
\numberwithin{equation}{section}

\begin{document}

\begin{center}
\textbf{\large{Derivative Polynomials and Closed-Form\\ Higher Derivative Formulae}}

\vskip 2mm\textbf{Djurdje Cvijovi\'{c}}\vskip 2mm

{\it Atomic Physics Laboratory, Vin\v{c}a Institute of
Nuclear Sciences \\
P.O. Box $522,$ $11001$ Belgrade$,$ Republic of Serbia}\\
\textbf{E-Mail: djurdje@vinca.rs}\\

\vskip 2mm \begin{quotation} \textbf{Abstract.} {\small In a recent paper,  Adamchik [V.S. Adamchik, On the Hurwitz function for rational arguments, {\it Appl. Math. Comp.}, {\bf 187} (2007), 3--12] expressed in a closed form symbolic derivatives of four functions  belonging to the class of functions whose derivatives are polynomials in terms of the same functions. In this sequel,  simple closed-form higher derivative formulae which involve the Carlitz-Scoville higher order tangent and secant numbers are derived for eight trigonometric and hyperbolic functions.}
\end{quotation}
\end{center}

\vskip 2mm\section{Introduction}
Recently, Adamchik \cite[Eqs. 26, 30, 31 and 32]{Adamchik} expressed  in a closed form symbolic derivatives of four functions  belonging to the class of functions whose derivatives are polynomials in terms of the same functions. In particular, he completely solved a long-standing problem of finding a closed-form expression for the higher derivatives of the cotangent function (see, for instance, \cite{Knuth,Apostol,Hoffman, Kolbig} and \cite[p. 161]{Berndt}) by showing that
\[\frac{\textup{d}^{n}
}{\textup{d}x^{n}}\,\cot(x)= (2\,\imath)^n \big(\cot(x)-\imath\big)\,\sum_{k\,=1}^{n}\frac{k!}{2}
{n\brace k} \,\big(\imath\,\cot(x)-1\big)^k,
\]
\noindent where $ {n \brace k}$ are the Stirling subset numbers.

In this sequel to the work of Adamchik, by using the derivative polynomials introduced by Hoffmann \cite{Hoffman,Hoffman2} (but see also references \cite{Knuth, Krishanamachary}), we further investigate the aforementioned class of functions and derive simple explicit closed-form higher derivative formulae for cotangent, tangent, cosecant and secant functions and their hyperbolic analogs. The formulae obtained involve the Carlitz-Scoville higher order tangent and secant numbers \cite{Carlitz,Carlitz2}.

\vskip 2mm\section{Derivative Polynomials}

Following Hoffman \cite{Hoffman,Hoffman2} we, by means of the exponential generating functions, define two sequences of polynomials, $\{P_n(x)\}_{n=0}^{\infty}$ and $\{Q_n(x)\}_{n=0}^{\infty}$, $ n\in \mathbb{N}_{0}:=\mathbb{N}\cup\{0\}, \textup{where}\,\, \mathbb{N}:=\{1, 2, 3,
\ldots\}$, which are here referred to as the derivative polynomials for tangent
\begin{equation}
P(x,t):=\frac{x+\tan (t)}{1-x \tan (t)}= \sum_{n=0}^{\infty}
P_n(x)\,\frac{t^n}{n!}
\end{equation}
\noindent and the derivative polynomials  for secant
\begin{equation}
Q(x,t):=\frac{\sec (t)}{1-x \tan (t)}= \sum_{n=0}^{\infty}
Q_n(x)\,\frac{t^n}{n!}.
\end{equation}
Equivalently, they may be defined by the formulae
\begin{equation}
\frac{\textup{d}^{n}
}{\textup{d}x^{n}}\,\tan(x)= P_n\big(\tan(x)\big)\tag{2.1*}
\end{equation}
\noindent and
\begin{equation}
\frac{\textup{d}^{n}
}{\textup{d}x^{n}}\,\sec(x)= \sec(x)\,Q_n\big(\tan(x)\big),\tag{2.2*}
\end{equation}
\noindent as it can be easily shown that
\begin{align*}
P_n\big(\tan(x)\big)& =\left.\frac{\textup{d}^{n}
}{\textup{d}t^{n}}\, P\big(\tan(x),t\big)\right|_{t = 0} = \left. \frac{\textup{d}^{n}
}{\textup{d}t^{n}}\tan\left(x+ t\right)\,\right|_{t = 0}
\\
&=\left. \frac{\textup{d}^{n}
}{\textup{d}x^{n}}\tan\left(x+ t\right)\,\right|_{t = 0}=\frac{\textup{d}^{n}
}{\textup{d}t^{n}}\tan(x)
\end{align*}
\noindent and
\begin{align*}
& \sec(x)\,Q_n\big(\tan(x)\big)= \sec(x)\,\left. \frac{\textup{d}^{n}
}{\textup{d}t^{n}}\,\frac{ \sec\left[\arctan\big(\tan(x)\big)+t\right]}{\sec\left[\arctan\big(\tan(x)\big)\right]} \,\right|_{t = 0}
\\
&=\left. \frac{\textup{d}^{n}
}{\textup{d}t^{n}}\,\sec\left(x+t\right)\,\right|_{t = 0}= \left. \frac{\textup{d}^{n}
}{\textup{d}x^{n}}\,\sec\left(x+t\right)\,\right|_{t = 0} =\frac{\textup{d}^{n}
}{\textup{d}x^{n}}\sec(x).
\end{align*}

\medskip
By making use of the chain rule it follows from (2.1*) that $P_{n}(x)$ satisfy
\begin{align}
&P_{0}(x) = x,\qquad P_{n}(x)=(1+x^2)\, P_{n-1}^{'}(x)\qquad(n\in\mathbb{N})\tag{2.1**}
\\ \nonumber\medskip
&\hspace{-7mm}\textup{and, similarly, from (2.2*) that}\nonumber
\\
&Q_{0}(x) = 1,\qquad Q_{n}(x)=(1+x^2)\, Q_{n-1}^{'}(x)+ x \,Q_{n-1}(x)\qquad(n\in\mathbb{N}).\tag{2.2**}
\end{align}
Another important and readily deducible  property of $P_{n}(x)$ and $Q_{n}(x)$ is
\begin{equation} P_{n}(-x)= (-1)^{n+1}\,P_{n}(x)\quad\textup{and}\quad Q_{n}(-x)= (-1)^{n}Q_{n}(x)\quad(n\in \mathbb{N}_{0}).
\end{equation}
\medskip
Upon noting that $\tan(x+\tfrac{\pi}{2})=-\cot(x)$ and $\sec(x+\tfrac{\pi}{2})=-\csc(x)$ and using (2.1*) and (2.2*) in conjunction with (2.5), we obtain
\begin{equation}
\frac{\textup{d}^{n}
}{\textup{d}x^{n}}\,\cot(x)= (-1)^n\,P_n\big(\cot(x)\big)
\end{equation}
\noindent and
\begin{equation}
\frac{\textup{d}^{n}
}{\textup{d}x^{n}}\,\csc(x)= (-1)^n\,\csc(x)\,Q_n\big(\cot(x)\big).
\end{equation}

We also consider the hyperbolic analogs of the derivative polynomials and define the derivative polynomials for hyperbolic tangent
\begin{equation}
\boldsymbol{\mathcal{P}}(x,t):=\frac{x+\tanh (t)}{1+x \tanh (t)}=\sum_{n=0}^{\infty}
\boldsymbol{\mathcal{P}}_n(x)\,\frac{t^n}{n!}
\end{equation}
\noindent and  for hyperbolic secant
\begin{equation}
\boldsymbol{\mathcal{Q}}(x,t):=\frac{\textup{sech}(t)}{1+x \tanh (t)}=\sum_{n=0}^{\infty}
\boldsymbol{\mathcal{Q}}_n(x)\,\frac{t^n}{n!},
\end{equation}
\noindent or, alternatively, as follows
\begin{equation}
\frac{\textup{d}^{n}
}{\textup{d}x^{n}}\,\tanh(x)= \boldsymbol{\mathcal{P}}_n\big(\tanh(x)\big)\tag{2.6*}
\end{equation}
\noindent and
\begin{equation}
\frac{\textup{d}^{n}
}{\textup{d}x^{n}}\,\textup{sech}(x)= \textup{sech}(x)\,\boldsymbol{\mathcal{Q}}_n\big(\tanh(x)\big).\tag{2.7*}
\end{equation}
These polynomials also may be generated by recurrence relations
\begin{align}
&\boldsymbol{\mathcal{P}}_{0}(x) = x,\qquad \boldsymbol{\mathcal{P}}_{n}(x)= (1-x^2)\, \boldsymbol{\mathcal{P}}_{n-1}^{'}(x)\qquad(n\in\mathbb{N})\tag{2.6**}
\\ \nonumber\medskip
&\hspace{-7mm}\textup{and}\nonumber
\\
&\boldsymbol{\mathcal{Q}}_{0}(x) = 1,\qquad \boldsymbol{\mathcal{Q}}_{n}(x)=(1-x^2)\, \boldsymbol{\mathcal{Q}}_{n-1}^{'}(x)- x \,\boldsymbol{\mathcal{Q}}_{n-1}(x)\qquad(n\in\mathbb{N})\tag{2.7**}
\end{align}
\noindent and they satisfy the symmetry relations
\begin{equation} \boldsymbol{\mathcal{P}}_{n}(-x)= (-1)^{n+1}\,\boldsymbol{\mathcal{P}}_{n}(x)\quad\textup{and}\quad \boldsymbol{\mathcal{Q}}_{n}(-x)= (-1)^{n}\boldsymbol{\mathcal{Q}}_{n}(x)\quad(n\in \mathbb{N}_{0}).
\end{equation}
Note that
\begin{equation}
\frac{\textup{d}^{n}
}{\textup{d}x^{n}}\,\coth(x)= \boldsymbol{\mathcal{P}}_n\big(\coth(x)\big)\quad\textup{and}\quad\frac{\textup{d}^{n}
}{\textup{d}x^{n}}\,\textup{csch}(x)= \textup{csch}(x)\,\boldsymbol{\mathcal{Q}}_n\big(\coth(x)\big).
\end{equation}

\vskip 4mm\section{Higher Derivative Formulae}

The {\em tangent numbers (of order} $k$)  $T(n,k)$ and {\em secant numbers (of order} $k$)  $S(n,k)$ are respectively defined by (see \cite[p. 428]{Carlitz} and \cite[p. 305]{Carlitz2})
\begin{equation}
\tan^{k}(t)=\sum_{n = k}^{\infty}
T(n,k)\,\frac{t^n}{n!}\qquad(k\in\mathbb{N})
\end{equation}
\noindent and
\begin{equation}
\sec(t) \tan^{k}(t)=\sum_{n = k}^{\infty}
S(n,k)\,\frac{t^n}{n!}\qquad(k\in\mathbb{N}_{0}).
\end{equation}
\noindent It is obvious, by parity considerations,  that $T(n,k)\neq 0$ is only when  $1\leq k\leq n $ and either both $n$ and $k$ are even or both $n$ and $k$ are odd. The same applies to $S(n,k)$ when $0\leq k\leq n $. Observe that,  $T(n,1)$ and $S(n,0)$ are, in fact, well-known the tangent and Euler numbers. Moreover, by (2.1) and (3.1), we have
\begin{equation}
P_n(0)= T(n,1)\qquad\textup{and}\qquad Q_n(0)= S(n,0).
\end{equation}

\medskip
Our main results are as follows.

\vskip 2mm \noindent{\bf Theorem 1.} {\em  Assume that } $n$ {\em and} $k$ {\em are nonnegative integers and let} $P_n(x)$ {\em and} $Q_n(x)$ {\em be the polynomials as defined by} \textup{(2.1)} {\em and} \textup{(2.2)}.
{\em Then, in terms of the tangent numbers of order} $k,$ $T(n,k),$ {\em given by} \textup{(3.1),}   {\em we have:}
\begin{equation}
P_n(x)= T(n,1) +\sum_{k\,=1}^{n+1} \frac{1}{k}\,T(n+1,k)\,x^k,
\end{equation}
\noindent {\em and, in terms of the secant numbers of order} $k,$  $S(n,k),$ {\em given by} \textup{(3.2),}   {\em we have:}
\begin{equation}
Q_n(x)=\sum_{k\,=0}^{n} S(n,k)\,x^k.
\end{equation}

\vskip 2mm \noindent {\bf Theorem 2.} {\em Let} $\boldsymbol{\mathcal{P}}_n(x)$ {\em and} $\boldsymbol{\mathcal{Q}}_n(x)$ {\em be the polynomials defined by} \textup{(2.6)} {\em and} \textup{(2.7)}. {\em Then:}
\begin{equation}
\boldsymbol{\mathcal{P}}_n(x)= (-1)^{\frac{n-1}{2}}\,\,T(n,1) +\sum_{k\,=1}^{n+1} \frac{(-1)^{\frac{n+k-1}{2}}}{k}\,\,T(n+1,k)\,x^k
\end{equation}
\noindent {\em and}
\begin{equation}
\boldsymbol{\mathcal{Q}}_n(x)=\sum_{k\,=0}^{n} (-1)^{\frac{n+k}{2}}\, S(n,k)\,x^k,
\end{equation}
\noindent {\em where} $T(n,k)$ {\em and} $S(n,k)$ {\em are the numbers} \textup{(3.1)} {\em and} \textup{(3.2)}.

\vskip 2mm \noindent{\bf Corollary 1.} {\em In terms of the tangent and  secant numbers of order} $k,$  $T(n,k)$ {\em and} $S(n,k),$ {\em for} $n\in \mathbb{N}_{0},$ {{\em we have:}
\begin{align*}
& \textup{(a)}\quad\frac{\textup{d}^{n}
}{\textup{d}x^{n}}\,\tan(x)= T(n,1)+\sum_{k\,=1}^{n+1}\frac{1}{k}\,T(n+1,k)\tan^k(x);
\\
&\textup{(b)}\quad\frac{\textup{d}^{n}
}{\textup{d}x^{n}}\,\sec(x)= \sec(x)\,\sum_{k\,=0}^{n} S(n,k)\,\tan^k(x);
\\
&\textup{(c)}\quad\frac{\textup{d}^{n}
}{\textup{d}x^{n}}\,\cot(x)= (-1)^n\,\left[T(n,1) +\sum_{k\,=1}^{n+1} \frac{1}{k}\,T(n+1,k)\,\cot^k(x)\right];
\\
&\textup{(d)}\quad\frac{\textup{d}^{n}
}{\textup{d}x^{n}}\,\csc(x)= (-1)^n\,\csc(x)\,\sum_{k\,=0}^{n} S(n,k)\,\cot^k(x).
\end{align*}

\vskip 2mm \noindent{\bf Corollary 2.} {\em In terms of the tangent and  secant numbers of order} $k,$ $T(n,k)$ {\em and} $S(n,k),$ {\em for} $n\in \mathbb{N}_{0},$ {{\em we have:}
\begin{align*}
& \textup{(a)}\quad\frac{\textup{d}^{n}
}{\textup{d}x^{n}}\tanh(x)=(-1)^{\frac{n-1}{2}}\,T(n,1)+\sum_{k\,=1}^{n+1} \frac{(-1)^{\frac{n+k-1}{2}}}{k}\,T(n+1,k)\tanh^k(x);
\\
&\textup{(b)}\quad\frac{\textup{d}^{n}
}{\textup{d}x^{n}}\,\textup{sech}(x)= \textup{sech}(x)\,\sum_{k\,=0}^{n}(-1)^{\frac{n+k}{2}}\, S(n,k)\,\tanh^k(x);
\\
&\textup{(c)}\quad\frac{\textup{d}^{n}
}{\textup{d}x^{n}}\,\coth(x)= (-1)^{\frac{n-1}{2}}\, T(n,1) +\sum_{k\,=1}^{n+1} \frac{(-1)^{\frac{n+k-1}{2}}}{k}\,T(n+1,k)\,\coth^k(x);
\\
&\textup{(d)}\quad\frac{\textup{d}^{n}
}{\textup{d}x^{n}}\,\textup{csch}(x)= \textup{csch}(x)\,\sum_{k\,=0}^{n}(-1)^{\frac{n+k}{2}}\, S(n,k)\,\coth^k(x).
\end{align*}

\vskip 2mm \noindent{\bf Remark 1.} We remark that, since, as detailed above, $T(n,k)$ and $S(n,k)$ are nonzero only under certain conditions, then the above given formulae can be written  (for computational purposes) in somewhat simplified form. For instance, we have

\begin{equation*}
P_{2 m-1}(x)= T(2 m-1,1) +\sum_{r\,=1}^{m} \frac{1}{2\,r}\,T(2 m,2 r)\,x^{2 r}\qquad(m\in\mathbb{N})
\end{equation*}
\noindent and
\begin{equation*}
P_{2 m}(x)= \sum_{r\,=0}^{m} \frac{1}{2\,r+1}\,T(2 m+1,2 r+1)\,x^{2 r+1}\qquad(m\in \mathbb{N}_{0}).
\end{equation*}

\vskip 1 mm \noindent {\bf Proof of Theorem 1.} In order to prove the formula (3.4) we first note that the generating function of the polynomials $P_n(x)$ can be rewritten as
\begin{equation*}
P(x,t) = \big(x+\tan (t)\big)\,\sum_{k=0}^{\infty}
\big( x\,\tan(t)\big)^k=x+\big(1+x^2\big)\,\sum_{k=1}^{\infty}
x^{k-1} \tan^k(t)
\end{equation*}
\noindent which, by making use of the definition of $T(n,k)$ in (3.1) and the elementary double series identities \cite[p. 57, Eq. (2)]{Rainville}
\[\sum_{n\,=1}^{\infty}\sum_{k\,=1}^{n}A(k,n)=\sum_{n\,=1}^{\infty}\sum_{k\,=1}^{\infty}A(k,n+k)=
\sum_{n\,=1}^{\infty}\sum_{k\,=n}^{\infty}A(k,n),
\]
\noindent  becomes
\begin{align}
P(x,t) &= x + \big(1+x^2\big)\,\sum_{k=1}^{\infty}
x^{k-1} \,\sum_{n=k}^{\infty}
T(n,k)\,\frac{t^n}{n!}\nonumber
\\
&= x + \sum_{n=1}^{\infty}\big(1+x^2\big)\,\left(\sum_{k=1}^{n}
T(n,k)\,x^{k-1}\right)\,\frac{t^n}{n!}.
\end{align}
On the other hand, by  (2.1) in conjunction with the recurrence (2.1**), we have
\begin{equation} P(x,t)= P_{0}(x)+\sum_{n=1}^{\infty} \,P_{n}(x)\,\frac{t^n}{n!}= P_{0}(x)+ \sum_{n=1}^{\infty} \big(1+x^2\big)\,P_{n-1}^{'}(x)\,\frac{t^n}{n!}
\end{equation}
\noindent and thus comparing (3.8) with  (3.9) clearly yields
\begin{equation*}P_{n-1}^{'}(x)= \sum_{k=1}^{n}
T(n,k)\,x^{k-1}
\end{equation*}
\noindent so that we find by integration that
\begin{equation}P_{n}(x)= P_{n}(0)+ \sum_{k=1}^{n+1}
\frac{1}{k}\,T(n+1,k)\,x^{k}.
\end{equation}
\noindent Now, in view of (3.3), the desired result (3.4) follows from (3.10).

Similarly, along the same lines, we have
\begin{align*}
Q(x,t) &= \sum_{k=0}^{\infty} \sec(t) \tan^k(t)\, x^k = \sum_{k=0}^{\infty} \left(\sum_{n=k}^{\infty}
S(n,k)\,\frac{t^n}{n!}\right) \,x^k
\\
&=\sum_{n=0}^{\infty}\left(\sum_{k=0}^{n} S(n,k)\,x^k \right) \frac{t^n}{n!}= \sum_{n=0}^{\infty} Q_n(x) \, \frac{t^n}{n!},
\end{align*}
\noindent and in this way we arrive at the second needed result (3.5).

\vskip 2 mm \noindent {\bf Proof of Theorem 2.} We first verify that
\[
\boldsymbol{\mathcal{P}}_n(x)=\imath^{n-1} P_n(\imath\,x)\qquad\textup{and}\qquad\boldsymbol{\mathcal{Q}}_n(x)=\imath^{n} Q_n(\imath\,x),
\]
\noindent and then, upon applying these last identities and Theorem 1, the two assertions of Theorem 2 follow.

\vskip 2 mm \noindent {\bf Proof of Corollaries 1 and 2.} The parts (a)--(d) of Corollary 1, in view of Theorem 1, are direct consequences of, respectively, the formulae (2.1*), (2.2*), (2.4) and (2.5). Similarly, the parts (a)--(d) of Corollary 2 follow by Theorem 2 and, respectively, formulae (2.6*), (2.7*) and (2.9).

\vskip 2mm
\section{Concluding Remarks}

We have explicitly expressed the trigonometric and hyperbolic derivative polynomials in a closed form in the terms of the higher tangent and secant numbers,  $T(n,k)$ and $S(n,k)$. The first few derivative polynomials are
\begin{align}
&\left.\begin{array}{l}P_{0}(x)\\\boldsymbol{\mathcal{P}}_{0}(x)\end{array}\right\}= x, \qquad \left.\begin{array}{l}P_{1}(x)\\\boldsymbol{\mathcal{P}}_{1}(x)\end{array}\right\}= \pm x^2+1,\qquad \left.\begin{array}{l}P_{2}(x)\\\boldsymbol{\mathcal{P}}_{2}(x)\end{array}\right\}= 2 x^3\pm 2 x,\hspace{10mm}\nonumber
\\
&\left.\begin{array}{l}P_{3}(x)\\\boldsymbol{\mathcal{P}}_{3}(x)\end{array}\right\}= \pm 6 x^4 + 8 x^2 \pm 2, \qquad \left.\begin{array}{l}P_{4}(x)\\\boldsymbol{\mathcal{P}}_{4}(x)\end{array}\right\}= 24 x^5 \pm 40 x^3 +16 x,\nonumber
\\
& \left.\begin{array}{l}P_{5}(x)\\\boldsymbol{\mathcal{P}}_{5}(x)\end{array}\right\}= \pm 120 x^6 + 240 x^4 \pm 136 x^2+16,\nonumber
\\
&\left.\begin{array}{l}P_{6}(x)\\\boldsymbol{\mathcal{P}}_{6}(x)\end{array}\right\}= 720 x^7 \pm 1680 x^5+ 1232 x^3\pm272 x, \nonumber
\end{align}
\noindent and
\begin{align}
&\left.\begin{array}{l}Q_{0}(x)\\\boldsymbol{\mathcal{Q}}_{0}(x)\end{array}\right\}=1, \qquad \left.\begin{array}{l}Q_{1}(x)\\\boldsymbol{\mathcal{Q}}_{1}(x)\end{array}\right\}= \pm x,\qquad \left.\begin{array}{l}Q_{2}(x)\\\boldsymbol{\mathcal{Q}}_{2}(x)\end{array}\right\}= 2 x^2 \pm 1,\hspace{20mm}\nonumber
\\
& \left.\begin{array}{l}Q_{3}(x)\\\boldsymbol{\mathcal{Q}}_{3}(x)\end{array}\right\}= \pm 6 x^3 + 5 x, \qquad \left.\begin{array}{l}Q_{4}(x)\\\boldsymbol{\mathcal{Q}}_{4}(x)\end{array}\right\}= 24 x^4 \pm 28 x^2 +5,\nonumber
\\
&\left.\begin{array}{l}Q_{5}(x)\\\boldsymbol{\mathcal{Q}}_{5}(x)\end{array}\right\}= \pm 120 x^5 + 180 x^3\pm61 x,\nonumber
\\
& \left.\begin{array}{l}Q_{6}(x)\\\boldsymbol{\mathcal{Q}}_{6}(x)\end{array}\right\}= 720 x^6 \pm 1320 x^4 + 662 x^2 \pm 61.\nonumber
\end{align}

It should be noted that the numbers $T(n,k)$ and $S(n,k)$ appear to be insufficiently investigated but simplicity of the the above-found formulae suggests that they would be interesting ones and well worthy of further study.

In conclusion, we confine ourselves to give only one of numerous consequences of the results presented in Section 3. The following formula
\begin{equation*}
\aligned
&\zeta\big(n,1-x\big) +(-1)^n\, \zeta\big(n,1-x\big)
\\
&=\frac{ (-1)^{n} \pi^n}{(n-1)!}\,\left[T(n-1,1) +\sum_{k\,=1}^{n} \frac{1}{k}\,T(n,k)\,\cot^k(\pi x)\right](n\in\mathbb{N}\setminus\{1\}; 0<x<1)
\endaligned
\end{equation*}
\noindent (\em c.f.} \cite[p. 8, Theorem 2.2]{Adamchik}) is obtained from our Corollary 1(c) and the reflection formula for the Hurwitz zeta function $\zeta(s,a)$ (see, for instance, \cite[Sec. 2.2]{Srivastava}).

\vskip 2mm

\noindent {\bf Remark.} Since submitting this paper (June 7, 2008) the author has learned of the following related publications: Refs. $\left[13–-15\right]$


\vskip 10mm
\begin{thebibliography}{20}

\bibitem {Adamchik} V.S. Adamchik, On the Hurwitz function for
rational arguments, {\it Appl. Math. Comp.}, {\bf 187} (2007), 3--12.

\bibitem{Knuth} D.E. Knuth and T.J. Buckholtz, Computation of tangent, Euler, and Bernoulli numbers, {\it Math. Comp.}, {\bf 21} (1967), 663–-688.

\bibitem {Apostol} T.M.  Apostol, Dirichlet L-functions and
character power sums, {\it J. Numb. Theory}, {\bf 2} (1970), 223--234.

\bibitem{Hoffman} M.E. Hoffman, Derivative polynomials for tangent and secant, {\it Amer. Math. Monthly}, {\bf 102} (1995), 23--30.

\bibitem{Kolbig} K.S. K\"{o}lbig, The polygamma function and the derivatives of the
cotangent function for rational arguments, CERN-IT-Reports,
CERN-CN-96-005, 1996.

\bibitem {Berndt} B.C.  Berndt, {\em Ramanujan's Notebooks}, Part I, Springer Verlag, New York, 1985.

\bibitem{Hoffman2} M.E.  Hoffman, Derivative polynomials, Euler polynomials, and associated integer sequences, {\it Electron. J. Combin.}, {\bf 6} (1999), R21.

\bibitem{Krishanamachary} C. Krishanamachary and R.M.  Bhimsenrao, On a table for calculating Eulerian numbers based on a new method, {\it Proc. London Math. Soc.}, Ser. 2, {\bf 22} (1923), 73--80.

\bibitem{Carlitz} L. Carlitz and R. Scoville, Tangent numbers and operators, {\it Duke Math. J.}, {\bf 39} (1972), 413--429.

\bibitem{Carlitz2} L. Carlitz, Permutations, sequences and special functions, {\it SIAM Review}, {\bf 17} (1975), 298--322.
\bibitem{Rainville} E.D. Rainville, {\it Special functions}, Macmillan, New York, 1960.

\bibitem{Srivastava} H.M.  Srivastava and J.  Choi, {\it Series Associated with the Zeta and Related Functions},
Kluwer Academic Publishers, Dordrecht, Boston and London, 2001.

\bibitem{Boyadzhiev} K.N. Boyadzhiev, Derivative polynomials for tanh, tan, sech and sec in explict form, Fibonacci Quart. 45 (2007),291–-303.

\bibitem{Franssens} G.R. Franssens, Functions with derivatives given by polynomials in the function itself or a related function, Analysis Mathematica 33 (2007), 17–-36.

\bibitem{Chang} C.-H. Chang, C.-W. Ha, Central factorial numbers and values of Bernoulli and Euler polynomials at rationals, Numer. Funct. Anal. Optimiz. 30 (2009), 214–-226.
\end{thebibliography}
\end{document}